                   \documentstyle[12pt,fleqn]{article}
\begin{document}\begin{flushright}\thispagestyle{empty}
OUT--4102--79\\
7 May 1999                                       \end{flushright}\vspace*{2mm}
                                                 \begin{center}{ \Large\bf
Parallel Integer Relation Detection:             \\[3pt]
Techniques and Applications$^{1)}$               }\vglue 10mm{\large\bf
David H.~Bailey$^{2)}$
and David J.~Broadhurst$^{3)}$                   }\end{center}\vfill
                                                 \noindent{\bf Abstract}\quad
Let $\{x_1, x_2, \cdots, x_n\}$ be a vector of real numbers.  An {\em
integer relation algorithm} is a computational scheme to find the $n$
integers $a_k$, if they exist, such that $a_1 x_1 + a_2 x_2 + \cdots +
a_n x_n= 0$.  In the past few years, integer relation algorithms have
been utilized to discover new results in mathematics and physics.
Existing programs for this purpose require very large amounts of
computer time, due in part to the requirement for multiprecision
arithmetic, yet are poorly suited for parallel processing.
This paper presents a new integer relation algorithm designed for
parallel computer systems, but as a bonus it also gives superior
results on single processor systems.  Single- and multi-level
implementations of this algorithm are described, together with
performance results on a parallel computer system.  Several
applications of these programs are discussed, including some new
results in number theory, quantum field theory and chaos theory.
\vfill\footnoterule\noindent{\small
$^1$) This work was supported by the Director, Office of Computational
and Technology Research, Division of Mathematical, Information, and
Computational Sciences of the U.S. Department of Energy, under contract
number DE-AC03-76SF00098.\\
$^2$) Lawrence Berkeley Laboratory, MS 50B-2239, Berkeley, CA 94720, USA\\
{\tt dhbailey@lbl.gov}\\
$^3)$ Open University, Department of Physics, Milton Keynes MK7 6AA, UK\\
{\tt D.Broadhurst@open.ac.uk}}
\newcommand{\nint}{{\rm nint }}
\newcommand{\li}{{\rm Li }}
\newcommand{\real}{{\rm Real }}
\newpage
\noindent {\bf 1. Introduction}

Let $x = (x_1, x_2, \cdots, x_n)$ be a vector of real numbers.  $x$ is
said to possess an integer relation if there exist integers $a_i$, not
all zero, such that $a_1 x_1 + a_2 x_2 + \cdots + a_n x_n = 0$.  By an
{\em integer relation algorithm}, we mean a practical computational
scheme that can recover (provided the computer implementation has
sufficient numeric precision) the vector of integers $a_i$, if it
exists, or can produce bounds within which no integer relation exists.

The problem of finding integer relations among a set of real numbers
was first studied by Euclid, who gave an iterative scheme which, when
applied to two real numbers, either terminates, yielding an exact
relation, or produces an infinite sequence of approximate relations.
The generalization of this problem for $n > 2$ was attempted by Euler,
Jacobi, Poincar\'e, Minkowski, Perron, Brun, Bernstein, among others.
The first integer relation algorithm with the required properties
mentioned above was discovered in 1977 by Ferguson and Forcade
\cite{ferfor}.  Since then, a number of other integer relation
algorithms have been discovered, including the ``HJLS'' algorithm
\cite{hjls} (which is based on the LLL algorithm), and the ``PSLQ''
algorithm.

\vspace{2ex} \noindent
{\bf 2. The PSLQ Algorithm}

The PSLQ integer relation algorithm features excellent numerical
stability, and it is effective in recovering a relation when the input
is known to only limited precision.  It has been generalized to
complex and even quaternion number systems.  A detailed discussion of
the PSLQ algorithm, together with a proof that the algorithm is
guaranteed to recover a relation in a polynomially bounded number of
iterations, is given in \cite{cpslq}.  The name ``PSLQ'' derives from
its usage of a partial sum of squares vector and a LQ
(lower-diagonal-orthogonal) matrix factorization.

A simple statement of the PSLQ algorithm, which is entirely equivalent
to the original formulation, is as follows: Let $x$ be the $n$-long
input real vector, and let nint denote the nearest integer function
(for exact half-integer values, define nint to be the integer with
greater absolute value).  Select $\gamma \geq \sqrt{4/3}$ (the authors
use $\gamma = \sqrt{4/3}$).  Then perform the following operations:

\vspace{2ex}

\noindent
Initialize:

\begin{enumerate}

\item For $j := 1$ to $n$: for $i := 1$ to $n$: if $i = j$ then set
$A_{i j} := 1$ and $B_{i j} := 1$ else set $A_{i j} := 0$ and $B_{i j}
:= 0$; endfor; endfor.

\item For $k := 1$ to $n$: set $s_k := \sqrt{\sum_{j=k}^n x_j^2}$;
endfor.  Set $t = 1 / s_1$.  For $k := 1$ to $n$: set $y_k := t x_k; \;
s_k := t s_k$; endfor.

\item Initial $H$: For $j := 1$ to $n - 1$: for $i := 1$ to $j - 1$:
set $H_{i j} := 0$; endfor; set $H_{j j} := s_{j+1} / s_j$; for $i :=
j + 1$ to $n$: set $H_{i j} := - y_i y_j / (s_j s_{j+1})$; endfor;
endfor.

\item Reduce $H$: For $i := 2$ to $n$: for $j := i - 1$ to $1$ step
$-1$: set $t := \nint (H_{i j} / H_{j j})$; and $y_j := y_j + t y_i$;
for $k := 1$ to $j$: set $H_{i k} := H_{i k} - t H_{j k}$; endfor; for
$k := 1$ to $n$: set $A_{i k} := A_{i k} - t A_{j k}$ and $B_{k j} :=
B_{k j} + t B_{k i}$; endfor; endfor; endfor.

\end{enumerate}

\noindent
Iteration: Repeat the following steps until precision has been
exhausted or a relation has been detected.

\begin{enumerate}

\item Select $m$ such that $\gamma^i |H_{i i}|$ is maximal when $i = m$.

\item Exchange the entries of $y$ indexed $m$ and $m + 1$, the
corresponding rows of $A$ and $H$, and the corresponding columns of
$B$.

\item Remove corner on $H$ diagonal: If $m \leq n - 2$ then set $t_0
:= \sqrt{H_{m m}^2 + H_{m,m+1}^2}$, $t_1 := H_{m m} / t_0$ and $t_2 :=
H_{m,m+1} / t_0$; for $i := m$ to $n$: set $t_3 := H_{i m}$, $t_4 :=
H_{i, m+1}$, $H_{i m} := t_1 t_3 + t_2 t_4$ and $H_{i, m+1} := - t_2
t_3 + t_1 t_4$; endfor; endif.

\item Reduce $H$: For $i := m + 1$ to $n$: for $j := \min(i - 1, m +
1)$ to $1$ step $-1$: set $t := \nint (H_{i j} / H_{j j})$ and $y_j :=
y_j + t y_i$; for $k := 1$ to $j$: set $H_{i k} := H_{i k} - t H_{j
k}$; endfor; for $k := 1$ to $n$: set $A_{i k} := A_{i k} - t A_{j k}$
and $B_{k j} := B_{k j} + t B_{k i}$; endfor; endfor; endfor.

\item Norm bound: Compute $M := 1 / \max_j |H_{j j}|$.  Then there can
exist no relation vector whose Euclidean norm is less than $M$.

\item Termination test: If the largest entry of $A$ exceeds the level
of numeric precision used, then precision is exhausted.  If the
smallest entry of the $y$ vector is less than the detection threshold
(see below), a relation has been detected and is given in the
corresponding column of $B$.

\end{enumerate}

It should be emphasized that for almost all applications of an integer
relation algorithm such as PSLQ, very high precision arithmetic must
be used.  Only a very small class of relations can be recovered
reliably with the 64-bit IEEE floating-point arithmetic that is
available on current computer systems.  In general, if one wishes to
recover a relation of length $n$, with coefficients of maximum size
$d$ digits, then it follows by an information theory argument that
the input vector $x$ must be specified to at least $n d$ digits, and
one must employ floating-point arithmetic accurate to at least $n d$
digits.  Practical integer relation programs always require greater
precision than this bound.  In fact, the difference between the level
of precision required for a given problem and the information theory
bound is a key figure of merit for integer relation algorithms.  PSLQ
is very efficient in this regard --- for most problems, PSLQ programs
can reliably recover relations with only about 15\% more digits of
precision than the information theory bound.

The software products Maple and Mathematica include multiple precision
arithmetic facilities.  One may also use any of several freeware
multiprecision software packages, such as the MPFUN package
(Fortran-77 and Fortran-90 versions are available), which was
developed by the first author \cite{mpfun,mpf90}, and the C/C++
version of MPFUN, which was recently developed by Sid Chatterjee and
Hermann Harjono of the University of North Carolina \cite{cmpfun}.
The two MPFUN packages permit one to write a program in conventional
Fortran-77/90 or C/C++, respectively, identifying some or all of the
variables to be multiple precision (integer, real or complex).  Then
in expressions where these variables appear, the appropriate multiple
precision routines are automatically referenced, thus saving
considerable programming effort.

In the course of the operation of the PSLQ algorithm on a real
computer system, the entries of the $y$ vector gradually decrease in
size, with the largest and smallest entries usually differing by no
more than two or three orders of magnitude.  When a relation is
detected by the algorithm, the smallest entry of the $y$ vector
abruptly decreases to roughly the multiprecision ``epsilon''
(i.e. $10^{-p}$, where $p > n d$ is the precision level in digits).
The detection threshold in the termination test (iteration step 6)
above is typically set to be a few orders of magnitude greater than
the epsilon value, in order to allow for reliable relation detection
in the presence of some numerical round-off error.  The ratio between
the smallest and the largest $y$ entry when a relation is detected can
be taken as a ``confidence level'' that the relation is a true
relation and not an artifact of insufficient numeric precision.  Very
small ratios at detection, such as $10^{-100}$, almost certainly
denote a true relation.

As shown in \cite{cpslq}, the PSLQ algorithm is guaranteed to find
relations in a bounded number of iterations.  However, this result is
based on the assumption of perfect, infinite-precision arithmetic.  In
an implementation on a real computer system, one can never rule out
hardware, software and programming errors, although the chances of
these errors can be minimized by independent computations.  Also, PSLQ
programs utilize multiprecision software with finite working
precision, and they make decisions based on numerical tolerances.
Thus it is possible that numerical anomalies can result, although
these anomalies generally can be remedied by using higher precision.

\vspace{2ex} \noindent
{\bf 3. Some Applications of the PSLQ Algorithm}

One application of PSLQ in the field of mathematical number theory is
to determine whether or not a given constant $\alpha$, whose value can
be computed to high precision, is algebraic of some degree $n$ or
less.  This can be done by first computing the vector $x = (1, \alpha,
\alpha^2, \cdots, \alpha^n$) to high precision and then applying an
integer relation algorithm.  If a relation is found for $x$, then this
relation vector is precisely the set of integer coefficients of a
polynomial satisfied by $\alpha$.  If a relation is not found, the
maximum bound determined by PSLQ means that $\alpha$ cannot be the
root of an polynomial of degree less than or equal to $n$, with
integer coefficients whose size (Euclidean norm) is less than the
established bound.  For example, it is well known \cite{pi-agm} that
\begin{eqnarray*}
\zeta(2) &=& 3 \sum_{k=1}^{\infty} \frac{1} {k^2 {2k \choose k}} \\
\zeta(3) &=& \frac{5}{2} \sum_{k=1}^{\infty} \frac{(-1)^{k-1}}
   {k^3 {2k \choose k}} \\
\zeta(4) &=& \frac{36}{17} \sum_{k=1}^{\infty} \frac{1}
   {k^4 {2k \choose k}}
\end{eqnarray*}

\noindent
These results have led some to suggest that
\begin{eqnarray*}
Z_5 &=& \zeta(5) / \sum_{k=1}^{\infty} \frac{(-1)^{k-1}}
   {k^5 {2k \choose k}}
\end{eqnarray*}

\noindent
might also be a simple rational or algebraic number.  Computations
using the PSLQ algorithm \cite{mpfun} have established that if $Z_5$
satisfies a polynomial of degree 25 or less, then the Euclidean norm
of the coefficients must exceed $2 \times 10^{37}$.  Results such as
this strongly suggest that the constants $\zeta(n)$ for $n > 4$ are
not given by simple one-term formulas as above.  Indeed, this
``negative'' result was fruitful in that it led to the discovery of
multi-term identities for such sums \cite{BorwBroad}.  An example will
be given in section 8.

One of the first ``positive'' results of this sort was the
identification of the constant $B_3 = 3.54409035955 \cdots$
\cite{mpfun}.  $B_3$ is the third bifurcation point of the logistic map
$x_{k+1} = r x_k (1 - x_k)$, which exhibits period doubling shortly
before the onset of chaos.  To be precise, $B_3$ is the smallest value
of the parameter $r$ such that successive iterates $x_k$ exhibit
eight-way periodicity instead of four-way periodicity.  Computations
using a predecessor algorithm to PSLQ found that $B_3$ is a root the
polynomial $4913 + 2108 t^2 - 604 t^3 - 977 t^4 + 8 t^5 + 44 t^6 + 392
t^7 - 193 t^8 - 40 t^9 + 48 t^{10} - 12 t^{11} + t^{12}$.  A stronger
result will be given in section 8.

A large number of results were recently found using PSLQ in the course
of research on multiple sums, such as those evaluated in Table 1.
After computing the numerical values of many of these constants, a
PSLQ program was used to determine if a given constant satisfied an
identity of a conjectured form.  These efforts produced numerous
empirical evaluations and suggested general results \cite{eulsum}.
Eventually, elegant proofs were found for many of these specific and
general results \cite{intrigue} and \cite{explicit}.  Three examples
of identities that are now proven are given in Table 1.  In the table,
$\zeta (t) = \sum_{j=1}^\infty j^{-t}$ is the Riemann zeta function,
and $\li_n (x) = \sum_{j=1}^\infty x^j j^{-n}$ denotes the
polylogarithm function.

\begin{table}[h]
$$
\begin{array}{|l|}
\hline
\sum_{k=1}^\infty \left(1 + \frac{1}{2} + \cdots + \frac{1}{k}\right)^2
  \, (k + 1)^{-4} \,=\, {37 \over 22680}\pi^6 - \zeta^2 (3) \\[3pt]
\sum_{k=1}^\infty \left(1 + \frac{1}{2} + \cdots + \frac{1}{k}\right)^3
  \, (k + 1)^{-6} \,=\,
\zeta^3(3)
 +\frac{197}{24} \zeta(9)
 +\frac{1}{2}\pi^2 \zeta(7)
    \\[3pt] \hspace{70mm}{}
 -\frac{11}{120} \pi^4 \zeta(5)
 -\frac{37}{7560} \pi^6 \zeta(3) \\[3pt]
\sum_{k=1}^\infty \left(1 - \frac{1}{2} + \cdots + (-1)^{k+1}
  \frac{1}{k}\right)^2 \, (k + 1)^{-3}
   \,=\, 4 \, \li_5(\frac12) - {1 \over 30} \ln^5(2) - {17 \over 32} \zeta(5)
    \\[3pt] \hspace{40mm}{}
   - {11 \over 720} \pi^4 \ln(2)
      + {7 \over 4} \zeta(3) \ln^2(2) + {1 \over 18} \pi^2 \ln^3(2)
      - {1 \over 8} \pi^2 \zeta(3) \\[3pt]\hline
\end{array}$$
\caption{Specimen evaluations, found with PSLQ and now proven}
\end{table}

It has been found that there is an intimate connection between such
multiple sums and the constants resulting from evaluation of Feynman
diagrams in quantum field theory \cite{BGK,BroadKreim}.  In
particular, the renormalization procedure (which removes infinities
from the perturbation expansion) entails {\em multiple zeta values\/}
defined by \cite{bbb}
\begin{eqnarray*}
\zeta (s_1, s_2, \cdots, s_r) &=& \sum_{k_1 > k_2 > \cdots > k_r > 0}
     \frac{1}{k_1^{s_1} \, k_2^{s_2} \cdots k_r^{s_r}}
\end{eqnarray*}
The $\zeta$ notation is used in analogy with Riemann's zeta function.
The PSLQ algorithm was used to find formulas and identities involving
these constants.  Again, a fruitful theory emerged, including a large
number of both specific and general results \cite{bbb,bbblc}.

Some recent quantum field theory results using PSLQ are even more
remarkable.  For example, it has now been shown \cite{broad98} that in
each of ten cases with unit or zero mass, the finite part the scalar
3-loop tetrahedral vacuum Feynman diagram reduces to 4-letter
``words'' that represent iterated integrals in an alphabet of 7
``letters'' comprising the one-forms $\Omega:=dx/x$ and
$\omega_k:=dx/(\lambda^{-k}-x)$, where $\lambda:=(1+\sqrt{-3})/2$ is
the primitive sixth root of unity, and $k$ runs from 0 to 5.  A
4-letter word is a 4-dimensional iterated integral, such as

\newpage

\begin{eqnarray*}
U &:=& \zeta(\Omega^2\omega_3\omega_0)=
    \int_0^1\frac{dx_1}{x_1}
    \int_0^{x_1}\frac{dx_2}{x_2}
    \int_0^{x_2}\frac{dx_3}{(-1-x_3)}
    \int_0^{x_3}\frac{dx_4}{(1-x_4)}\\
  &=&\sum_{j>k>0}\frac{(-1)^{j+k}}{j^3k}
\end{eqnarray*}
There are $7^4$ four-letter words.  Only two of these are primitive
terms occurring in the 3-loop Feynman diagrams: $U$, above, and
\begin{eqnarray*}
V &:=& \real[\zeta(\Omega^2\omega_3\omega_1)] = \sum_{j > k > 0}\frac{ (-1)^j
\cos (2 \pi k / 3) }{ j^3 k}.
\end{eqnarray*}
The remaining terms in the diagrams reduce to products of constants
found in Feynman diagrams with fewer loops.  These ten cases as shown
in Figure 1.  In these diagrams, dots indicate particles with nonzero
rest mass.  The formulas that have been found for the corresponding
constants are given in Table 2.  The constant $C = \sum_{k > 0} \sin
(\pi k / 3) / k^2$.

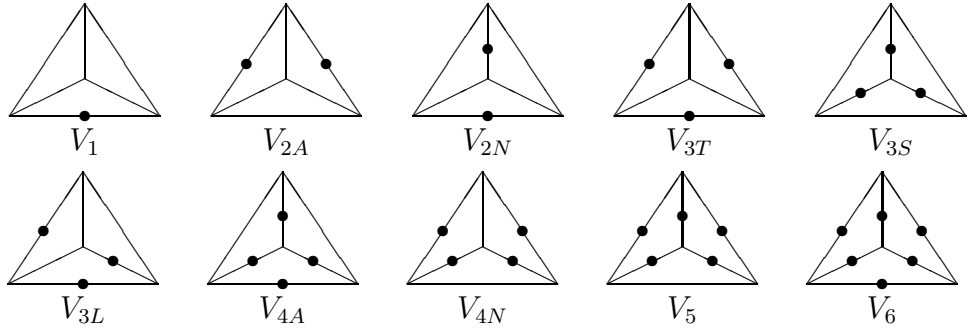
\begin{figure}[h]

\setlength{\unitlength}{0.01cm}
\newbox\shell
\newcommand{\dia}[2]{\setbox\shell=\hbox{\begin{picture}(180,220)(-90,-110)#1
\put(-90,-95){\makebox(180,220)[b]{$#2$}}\end{picture}}\dimen0=\ht
\shell\multiply\dimen0by7\divide\dimen0by16\raise-\dimen0\box\shell\hfill}
\newcommand{\mass}{\circle*{15}}
\newcommand{\blone}{\put(-53,20){\mass}}
\newcommand{\bltwo}{\put(-40,-20){\mass}}
\newcommand{\blthree}{\put(0,40){\mass}}
\newcommand{\blfour}{\put(53,20){\mass}}
\newcommand{\blfive}{\put(40,-20){\mass}}
\newcommand{\blsix}{\put(0,-50){\mass}}
\newcommand{\mrk}[3]{\put(#1,#2){\mbox{#3}}}
\newcommand{\tet}{\put(-100,-50){\line(1,0){200}}
\newcommand{\deq}{\put(200,0){\mbox{$=$}}}
\put(-100,-50){\line(2,3){100}}\put(100,-50){\line(-2,3){100}}
\put(-100,-50){\line(2,1){100}}\put(100,-50){\line(-2,1){100}}
\put(0,0){\line(0,1){100}}}
\newcommand{\tetmrk}[6]{\tet
\mrk{-75}{20}{#1}\mrk{-60}{-20}{#2}\mrk{0}{20}{#3}
\mrk{55}{20}{#4}\mrk{40}{-20}{#5}\mrk{-10}{-80}{#6}}

\mbox{\hspace{1cm}}\hfill
\dia{\tet\blsix}{V_1}
\dia{\tet\blone\blfour}{V_{2A}}
\dia{\tet\blthree\blsix}{V_{2N}}
\dia{\tet\blone\blfour\blsix}{V_{3T}}
\dia{\tet\bltwo\blthree\blfive}{V_{3S}}\mbox{\hspace{1cm}}\\
\mbox{\hspace{1cm}}\hfill
\dia{\tet\blone\blfive\blsix}{V_{3L}}
\dia{\tet\bltwo\blthree\blfive\blsix}{V_{4A}}
\dia{\tet\blone\bltwo\blfour\blfive}{V_{4N}}
\dia{\tet\blone\bltwo\blthree\blfour\blfive}{V_5}
\dia{\tet\blone\bltwo\blthree\blfour\blfive\blsix}{V_6}
\mbox{\hspace{1cm}}\par

\caption{The ten tetrahedral cases}
\end{figure}

\begin{table}[h]
$$\begin{array}{|lcl|}\hline
V_1 &=& 6 \zeta(3) + 3 \zeta(4) \\[3pt]
V_{2A} &=& 6 \zeta(3) - 5 \zeta(4) \\[3pt]
V_{2N} &=& 6 \zeta(3) - \frac{13}{2} \zeta(4) - 8 U \\[3pt]
V_{3T} &=& 6 \zeta(3) - 9 \zeta(4) \\[3pt]
V_{3S} &=& 6 \zeta(3) - \frac{11}{2} \zeta(4) - 4 C^2 \\[3pt]
V_{3L} &=& 6 \zeta(3) - \frac{15}{4} \zeta(4) - 6 C^2 \\[3pt]
V_{4A} &=& 6 \zeta(3) - \frac{77}{12} \zeta(4) - 6 C^2 \\[3pt]
V_{4N} &=& 6 \zeta(3) - 14 \zeta(4) - 16 U \\[3pt]
V_5 &=& 6 \zeta(3) - \frac{469}{27} \zeta(4) + \frac{8}{3} C^2 - 16 V \\[3pt]
V_6 &=& 6 \zeta(3) - 13 \zeta(4) - 8 U - 4 C^2\\[3pt]\hline
\end{array}$$
\caption{Evaluations of the 10 constants corresponding to the 10
cases in Figure 1}
\end{table}

\newpage

\vspace{2ex} \noindent
{\bf 4. A New Formula for Pi}

Through the centuries mathematicians have assumed that there is no
shortcut to computing just the $n$-th digit of $\pi$.  Thus, it came
as no small surprise when such an algorithm was recently discovered
\cite{bbp}.  In particular, this simple scheme allows one to compute
the $n$-th hexadecimal (or binary) digit of $\pi$ without computing
any of the first $n-1$ digits, without using multiple-precision
arithmetic software, and at the expense of very little computer
memory.  The one millionth hex digit of $\pi$ can be computed in this
manner on a current-generation personal computer in only about 60
seconds run time.

This scheme is based on the following new formula, which was
discovered using PSLQ:
\begin{eqnarray*}
\pi &=& \sum_{k=0}^\infty \frac{1}{16^k} \left[
 \frac{4}{8k+1} - \frac{2}{8k+4} - \frac{1}{8k+5} - \frac{1}{8k+6} \right]
\end{eqnarray*}

\noindent
It is likely the first instance in history of a significant new
formula for $\pi$ discovered by computer.  Further base-2 results
are given in \cite{bbp,broad-z3z5}. In \cite{broad98} base-3 results
were obtained, including
\begin{eqnarray*}\pi^2&=&\frac{2}{27}
\sum_{k=0}^\infty\frac{1}{729^k}\left[
 \frac{243}{(12k+1)^2}
-\frac{405}{(12k+2)^2}
-\frac{81}{(12k+4)^2}
-\frac{27}{(12k+5)^2}
\right.\nonumber\\&&\left.{}
-\frac{72}{(12k+6)^2}
-\frac{9}{(12k+7)^2}
-\frac{9}{(12k+8)^2}
-\frac{5}{(12k+10)^2}
+\frac{1}{(12k+11)^2}
\right]
\end{eqnarray*}

\vspace{2ex} \noindent
{\bf 5. Multi-Level Implementations of PSLQ}

In spite of the relative efficiency of PSLQ compared to the other
algorithms in the literature, computer run times of programs that
straightforwardly implement the PSLQ algorithm are typically quite
long.  Even modest-sized problems can require many hours for solution
on a current personal computer or workstation.  This is mainly due to
the cost of using high precision arithmetic software for nearly every
operation in the algorithm.

As it turns out, it is possible to perform most, if not all, of the
PSLQ iterations using ordinary 64-bit computer arithmetic, with only
occasional recourse to multiprecision arithmetic.  In this way, run
times can be dramatically reduced.  Here is a sketch of this scheme,
which will be referred to as a ``two-level'' implementation of the
PSLQ algorithm.  In the following, ``double precision'' means the
64-bit IEEE hardware arithmetic available on most current computer
systems, and $\bar{y}$, $\bar{A}$, $\bar{B}$ and $\bar{H}$ denotes double
precision counterparts to the arrays $y$, $A$, $B$ and $H$ in the PSLQ
algorithm.

First, perform the multiprecision initialization steps of PSLQ as
given in section 2 above.  Then perform a double precision
``re-initialization'' step: set $\bar{A}$ and $\bar{B}$ to the $n
\times n$ identity matrix; set $\bar{y}$ to the best double precision
approximation of the current $y$ vector, multiplied by a scale factor
so that its largest entry is unity; and set $\bar{H}$ to the best
double precision approximation of the current $H$ matrix.  For some
extremely large problems it may be necessary to scale the $\bar{H}$
matrix to avoid numeric overflow.  Then perform an LQ
(lower-diagonal-orthogonal) matrix factorization on $\bar{H}$, and
replace $\bar{H}$ by the lower diagonal portion of the result (the
upper right portion is zeroed).  The subroutine DQRDC of the Linpack
library \cite{linpack} may be employed for this factorization,
provided both the input and output matrices are transposed.

Next, perform PSLQ iterations using the double precision arrays.  In
the course of these iterations, the entries of $\bar{A}$ and $\bar{B}$
(which contain integer values, although stored as IEEE double
precision data), steadily increase in size.  Monitor the entries of
these matrices as they are updated, and when any entry reaches a
certain threshold (the authors use $10^{13}$), or when the smallest
$\bar{y}$ entry becomes smaller than a certain threshold (the authors
use $10^{-14}$), then update the multiprecision arrays by means of
matrix multiplication operations, as follows:
\begin{eqnarray*}
y &:=& y \cdot \bar{B} \\
B &:=& B \cdot \bar{B} \\
A &:=& \bar{A} \cdot A \\
H &:=& \bar{A} \cdot H
\end{eqnarray*}

\noindent
After these updates are performed, the entries of the $A$ matrix and
the $y$ vector are checked, as in the termination test (iteration step
6) of PSLQ, and a norm bound is computed.  If neither of the
termination conditions holds, then the double precision arrays are
re-initialized again as mentioned above, another set of double
precision iterations are performed, and the process continues.

This general scheme works well for many problems, but there are
several difficulties that must be dealt with in a fully robust
implementation.  One difficulty is that at some point in the
computation (typically at the very beginning), the $y$ vector may have
a dynamic range that exceeds the range (11 or 12 orders of magnitude)
that can be safely handled using double precision iterations.  Another
difficulty is that occasionally an entry is produced in the $\bar{A}$
or $\bar{B}$ matrix that exceeds the largest whole number ($2^{53} =
9.007 \cdots \times 10^{15})$ that can be exactly represented as
64-bit IEEE data.  A straightforward solution when such a condition
occurs is to abandon the current iteration, restore a previous
iteration's values of $\bar{y}$, $\bar{A}$, $\bar{B}$ and $\bar{H}$,
update the multiprecision arrays as above, perform an LQ matrix
factorization on the $H$ matrix, and then perform iterations using
full multiprecision arithmetic until these special conditions no
longer hold.

A more efficient solution for large problems that require very high
precision is to employ ``intermediate precision'', in other words a
fixed level of precision (the authors use 125 digits) that is
intermediate between double precision and full multiprecision.
Updating the full multiprecision arrays from the intermediate
precision arrays is done with matrix multiplication operations in a
manner precisely analogous to that described above.  Incorporating
intermediate precision in this manner gives rise to what we will refer
to as a ``three-level'' implementation of PSLQ.

One additional improvement that can be made to each of these schemes
is to omit multiprecision computation of the $A$ matrix (although the
double precision and intermediate precision equivalents of $A$ must be
computed).  The multiprecision $A$ matrix (which is the inverse of the
$B$ matrix) is used in the PSLQ algorithm only to determine when
execution must be halted due to the exhaustion of numeric precision.
However, exhaustion of numeric precision can alternatively be handled
by halting iterations when the smallest $y$ entry is sufficiently
close to the multiprecision epsilon level (the authors use a factor of
$10^{25}$).

These three PSLQ schemes (one-level, two-level and three-level) have
been implemented by the first author, using the Fortran-90 MPFUN
software \cite{mpf90}.  Some performance results are shown in Table 3
for a class of problems.  Here $r, s$ define the constant $\alpha =
3^{1/r} - 2^{1/s}$, which is algebraic of degree $rs$, and $n = rs +
1$.  The $n$-long vector of coefficients of the polynomial satisfied
by $\alpha$ can thus be obtained by using a PSLQ program, as explained
in section 3.  The column headed ``Iterations'' gives the number of
PSLQ iterations required for solution, while ``Digits'' gives the
working precision level used, in decimal digits.  ``Time'' gives CPU
time in seconds for runs on a single processor of an SGI Origin-2000
system with 195 MHz R10000 CPUs.

It can be seen from these results that the two-level PSLQ program is
up to 65 times faster than the one-level program, yet it finds
relations just as well, usually in exactly the same course of
iterations as the one-level program.  The three-level program is
faster than the two-level program for large problems, even though the
special conditions mentioned above rarely arise in the particular
problems mentioned in the table.  The reason for this fortunate
circumstance appears to be improved data locality in the three-level
scheme, which is advantageous on modern cache-based computer systems.
Fully detailed computer programs are available from the authors at the
web site {\tt http://www.nersc.gov/\~{}dhb} .

\begin{table} \begin{center}
\begin{tabular}{|l|r|r|r|r|r|r|r|r|}
\hline
& & & \multicolumn{2}{|c|}{One-level} & \multicolumn{2}{|c|}{Two-level} &
  \multicolumn{2}{|c|}{Three-level} \\
$r,s$ & $n$ & Iterations & Digits & Time & Digits & Time & Digits & Time \\
\hline
5,5   &  26 &   5143 & 180 &   32.37 &  190 &    1.29 &      & \\
5,6   &  31 &   9357 & 240 &  105.48 &  250 &    3.16 &      & \\
6,6   &  37 &  15217 & 310 &  298.85 &  320 &    7.19 &      & \\
6,7   &  43 &  25361 & 420 &  942.66 &  420 &   17.22 &      & \\
7,7   &  50 &  36947 & 500 & 2363.71 &  510 &   36.29 &      & \\
7,8   &  57 &  60817 &     &         &  680 &   90.08 &      & \\
8,8   &  65 &  86684 &     &         &  850 &  195.19 &  910 &  233.48 \\
8,9   &  73 & 124521 &     &         & 1050 &  425.67 & 1120 &  460.34 \\
9,9   &  82 & 174140 &     &         & 1310 &  934.96 & 1370 &  922.90 \\
9,10  &  91 & 245443 &     &         & 1620 & 2032.69 & 1680 & 1780.65 \\
10,10 & 101 & 342931 &     &         & 2000 & 4968.64 & 2060 & 3366.92 \\
\hline
\end{tabular} \end{center}
\caption{Run times for the three PSLQ programs}
\end{table}

\vspace{2ex} \noindent
{\bf 6. The Multi-Pair Algorithm}

Even with the substantial accelerations described in the previous
section, run times are painfully long for some very large problems of
current interest in mathematics and physics.  Thus one is led to
consider employing highly parallel supercomputers, which have the
potential of performance hundreds of times faster than for
single-processor scientific workstations and personal computers.

Unfortunately, the standard PSLQ algorithm appears singularly unsuited
for modern parallel computer systems, which require high levels of
coarse-grained concurrency.  The main difficulty is that large integer
relation problems often require over one million PSLQ iterations, each
of which must be completed before the next begins.  Further, within an
individual iteration, the key reduction operation (iteration step 4)
has a recursion that inhibits any possibility for parallel execution,
except at the innermost loop level.  These considerations have led
some researchers in the field to conclude that there is no hope for
any significant parallel acceleration of PSLQ-type computations.

But it turns out that a variant of the PSLQ algorithm can be
formulated that dramatically reduces the number of sequential
iterations that must be performed, while at the same time exhibiting
reasonably high concurrency in the major steps of individual
iterations.  To that end, consider the following algorithm, which will
be referred to as the ``multi-pair'' variant of PSLQ.  Here $\gamma =
\sqrt{4/3}$ as before, and $\beta = 0.4$.

\vspace{2ex}

\noindent
Initialize:

\begin{enumerate}

\item For $j := 1$ to $n$: for $i := 1$ to $n$: if $i = j$ then set
$A_{i j} := 1$ and $B_{i j} := 1$ else set $A_{i j} := 0$ and $B_{i j}
:= 0$; endfor; endfor.

\item For $k := 1$ to $n$: set $s_k := \sqrt{\sum_{j=k}^n x_j^2}$;
endfor; set $t = 1 / s_1$; for $k := 1$ to $n$: set $y_k := t x_k; \;
s_k := t s_k$; endfor.

\item Initial $H$: For $j := 1$ to $n - 1$: for $i := 1$ to $j - 1$:
set $H_{i j} := 0$; endfor; set $H_{j j} := s_{j+1} / s_j$; for $i :=
j + 1$ to $n$: set $H_{i j} := - y_i y_j / (s_j s_{j+1})$; endfor;
endfor.

\end{enumerate}

\noindent
Iteration: Repeat the following steps until precision has been
exhausted or a relation has been detected.

\begin{enumerate}

\item Sort the entries of the $(n-1)$-long vector $\{\gamma^i |H_{i
i}|\}$ in decreasing order, producing the sort indices.

\item Beginning at the sort index $m_1$ corresponding to the largest
$\gamma^i |H_{i i}|$, select pairs of indices $(m_i, m_i + 1)$,
where $m_i$ is the sort index.  If at any step either $m_i$ or $m_i +
1$ has already been selected, pass to the next index in the list.
Continue until either $\beta n$ pairs have been selected, or the list
is exhausted.  Let $p$ denote the number of pairs actually selected in
this manner.

\item For $i := 1$ to $p$, exchange the entries of $y$ indexed $m_i$
and $m_i + 1$, and the corresponding rows of $A$, $B$ and $H$; endfor.

\item Remove corners on $H$ diagonal: For $i := 1$ to $p$: if $m_i
\leq n - 2$ then set $t_0 := \sqrt{H_{m_i, m_i}^2 + H_{m_i, m_i +
1}^2}$, $t_1 := H_{m_i, m_i} / t_0$ and $t_2 := H_{m_i, m_i + 1} /
t_0$; for $i := m_i$ to $n$: set $t_3 := H_{i, m_i}$; $t_4 := H_{i,m_i
+ 1}$; $H_{i,m_i} := t_1 t_3 + t_2 t_4$; and $H_{i,m_i + 1} := - t_2
t_3 + t_1 t_4$; endfor; endif; endfor.

\item Reduce $H$: For $i := 2$ to $n$: for $j := 1$ to $n - i + 1$:
set $l := i + j - 1$; for $k := j + 1$ to $l - 1$: set $H_{l j} :=
H_{l j} - T_{l k} H_{k j}$; endfor; set $T_{l j} := \nint (H_{l j} /
H_{j j})$ and $H_{l j} := H_{l j} - T_{l j} H_{j j}$; endfor; endfor.

\item Update $y$: For $j := 1$ to $n - 1$: for $i := j + 1$ to $n$:
set $y_j := y_j + T_{i j} y_i$; endfor; endfor.

\item Update $A$ and $B$: For $k := 1$ to $n$: for $j := 1$ to $n -
1$: for $i := j + 1$ to $n$: set $A_{i k} := A_{i k} - T_{i j} A_{j
k}$ and $B_{j k} := B_{j k} + T_{i j} B_{i k}$; endfor; endfor;
endfor.

\item Norm bound: Compute $M := 1 / \max_j |H_{j j}|$.  Then there can
exist no relation vector whose Euclidean norm is less than $M$.

\item Termination test: If the largest entry of $A$ exceeds the level
of numeric precision used, then precision is exhausted.  If the
smallest entry of the $y$ vector is less than the detection threshold
(see section 2), a relation has been detected and is given in the
corresponding row of $B$.

\end{enumerate}

There are several differences between this algorithm and the standard
one-level PSLQ algorithm: (1) there is no reduction step in the
initialization; (2) the $B$ matrix is transposed from the standard
PSLQ algorithm; (3) up to $\beta n$ disjoint pairs (not just a single
pair) of adjacent indices are selected in each iteration; (4) the $H$
reduction loop proceeds along successive lower diagonals of the $H$
matrix; (5) a $T$ matrix is employed, which contains the $t$
multipliers of the standard PSLQ; and (6) the $y$, $A$ and $B$ arrays
are not updated with $H$, but in separate loops.

Since the multi-pair algorithm maintains the $H$ matrix in lower
triangular form, and the $A$ and $B$ matrices are maintained as
integer matrices, the norm bound stated in iteration step 8 above is
valid, by the same argument that applies to the original PSLQ
algorithm \cite{cpslq}.

Unfortunately, we cannot offer a proof that the multi-pair algorithm
is guaranteed to recover a relation in a bounded number of iterations,
as can be done with PSLQ.  In fact, it has been found that for certain
special problems, the multi-pair algorithm, as stated above, falls
into a repeating cycle, with a period of (usually) two iterations.
Our implementation deals with this difficulty by comparing the $y$
vector at the end of each iteration with saved copies from eight
previous iterations, and if a duplication is found, then only one pair
of indices is selected in step 2 of the next iteration (so that the
next iteration is equivalent to a standard PSLQ iteration).  It should
be added, however, that these repeating situations are extremely rare
in nontrivial problems.  We have not seen any instances of such
repeats when $n > 20$.

On the positive side, we have found, based on our experience with a
wide variety of sample problems, that the norm bound increases much
more rapidly than in the standard PSLQ.  Indeed, it appears that the
selection of up to $\beta n$ disjoint pairs of indices in step 2 above
has the effect of reducing the iteration count by nearly the factor
$\beta n$.  This results in a significant saving in the number of
expensive $H$ reduction and array update steps.  More importantly,
without this dramatic reduction in the sequential iteration count, an
efficient parallel implementation would not be possible.  Parallel
issues will be discussed in greater detail in the next section.

Given that the multi-level implementations of PSLQ are so much faster
than the standard one-level PSLQ, one might also wonder whether there
exist analogous multi-level implementations of the multi-pair
algorithm.  Happily, the multi-level scheme sketched in section 5 can
be adopted almost without change.  One change that is required is that
the multiprecision arrays are updated as follows:
\begin{eqnarray*}
y &:=& \bar{B} \cdot y \\
B &:=& \bar{B} \cdot B \\
A &:=& \bar{A} \cdot A \\
H &:=& \bar{A} \cdot H
\end{eqnarray*}

\noindent
Note that $y$ and $B$ are updated here in the same manner as the $A$
and $H$ arrays.  This change stems from the fact that the $B$ matrix
in the multi-pair scheme is transposed from the $B$ matrix in the
standard PSLQ algorithm.

The multi-pair algorithm and the multi-level implementations described
here were all devised to permit parallel processing.  But it turns out
that these programs also run faster on a single processor system,
compared with the standard PSLQ equivalents.  Some one-processor
timings are shown in Table 4 for the suite of test problems used in
Table 3.  Note for example that the one-level multi-pair program is up
to twice as fast as the one-level PSLQ program, and the three-level
multi-pair program is up to 22\% faster than the three-level PSLQ
program.  Note also that the iteration counts are reduced by a factor
of up to 34.  Finally, note that the multi-pair schemes require
slightly less numeric precision for solution than their PSLQ
counterparts.  The reason for this unanticipated benefit is not known.

\begin{table}[h]
\begin{center}
\begin{tabular}{|l|r|r|r|r|r|r|r|r|}
\hline
& & & \multicolumn{2}{|c|}{One-level} & \multicolumn{2}{|c|}{Two-level} &
  \multicolumn{2}{|c|}{Three-level} \\
$r,s$ & $n$ & Iterations & Digits & Time & Digits & Time & Digits & Time \\
\hline
5,5  &  26 &   558 & 180 &   26.08 &  180 &    1.48 & & \\
5,6  &  31 &   840 & 230 &   70.71 &  240 &    3.43 & & \\
6,6  &  37 &  1136 & 310 &  189.27 &  310 &    7.84 & & \\
6,7  &  43 &  1625 & 400 &  479.07 &  410 &   17.22 & & \\
7,7  &  50 &  2071 & 500 & 1130.85 &  500 &   35.64 & & \\
7,8  &  57 &  2410 &     &         &  660 &   69.39 & & \\
8,8  &  65 &  3723 &     &         &  800 &  169.62 &  880 &  214.66 \\
8,9  &  73 &  4943 &     &         & 1010 &  358.07 & 1100 &  427.29 \\
9,9  &  82 &  6169 &     &         & 1260 &  744.20 & 1320 &  804.51 \\
9,10 &  91 &  7850 &     &         & 1560 & 1556.37 & 1600 & 1450.29 \\
10,10& 101 & 10017 &     &         & 1890 & 3283.08 & 1950 & 2747.12 \\
\hline
\end{tabular} \end{center}
\caption{Run times for the three multi-pair programs}
\end{table}

\vspace{2ex} \noindent
{\bf 7. Parallel Implementations of the Multi-Pair Algorithm}

The key steps of the multi-pair iterations are all suitable for
parallel execution.  First note that the $p$ row exchanges in
iteration step 3, as well as the $p$ corner removal operations in step
4, can be performed concurrently, since the $p$ pairs of indices
$(m_i, m_i + 1)$ are all disjoint.  Secondly, the reorganized $H$
matrix reduction step (step 5), which is equivalent to the $H$ matrix
reduction scheme in the standard PSLQ, may be performed concurrently
at the second loop level, instead of only at the innermost loop level
as in standard PSLQ.  The update of the $A$ and $B$ arrays (step 7) is
even more favorable to parallel processing: this loop may be performed
concurrently at the outermost loop level.  The change in the $B$
matrix, which is transposed from the standard PSLQ algorithm, is
favorable for an implementation on a distributed memory parallel
computer.

The two- and three-level multi-pair schemes are also well suited for
parallel computation.  This is because the dominant cost of these
programs is the matrix multiplication operations involved in the
multiprecision array updates, and these matrix multiplications can be
performed concurrently at the outermost loop level.  The parallel
techniques mentioned in the previous paragraph can still be applied to
the double precision and intermediate precision iterations.  It turns
out, though, that the double precision iterations run so rapidly that
parallel processing of these iterations is often not worth the
overhead.  Nonetheless, we have achieved modest acceleration on very
large problems by using parallel processing on some steps of double
precision iterations.  Some parallel performance results will be given
in the next section.

\vspace{2ex} \noindent
{\bf 8. Large Applications and Parallel Performance}

Three recent applications will be described here, each of which
involves very large integer relation problems.  Thus they are
excellent test cases for the new multi-pair programs.

\vspace{1ex} \noindent {\bf Reduction of Euler sums:}
In section 3, we mentioned recent research on multiple zeta values,
which play a key role in quantum field theory \cite{BroadKreim}.  More
generally, one may define {\em Euler sums} by \cite{bbb}
\begin{eqnarray*}
\zeta\left(\begin{array}{cccc}
s_1,&s_2&\cdots&s_r\\
\sigma_1,&\sigma_2&\cdots&\sigma_r
\end{array}\right):=
\sum_{k_1 > k_2 > \cdots > k_r > 0}
     \frac{\sigma_1^{k_1}}{k_1^{s_1}}\,
     \frac{\sigma_2^{k_2}}{k_2^{s_2}}\,
     \cdots\,
     \frac{\sigma_r^{k_r}}{k_r^{s_r}}
\end{eqnarray*}
where $\sigma_j=\pm1$ are signs and $s_j>0$ are integers.  When all
the signs are positive, one has a multiple zeta value.  Constants with
alternating signs appear in problems such as computation of the
magnetic moment of the electron.

It was conjectured by the second author that the dimension of the
space of Euler sums with weight $w:=\sum_j s_j$ is the Fibonacci
number $F_{w+1} = F_w + F_{w-1}$, with $F_1 = F_2 = 1$.  Complete
reductions of all Euler sums to a basis of size $F_{w+1}$ have been
obtained with PSLQ at weights $w\leq9$.  At weights $w=10$ and $w=11$
the conjecture has been stringently tested by application of PSLQ in
more than 600 cases.  At weight $w=11$ such tests involve solving
integer relations of size $n = F_{12} + 1 = 145$.  In a typical case,
each of the 145 constants was computed to more than 5,000 digit
accuracy, and a working precision level of 5,000 digits was employed
in the three-level multi-pair program.  A relation was detected at
iteration 31,784.  The minimum and maximum $y$ vector entries at the
point of detection were $9.515 \times 10^{-4970}$ and $4.841 \times
10^{-4615}$, respectively.  The ratio of these two values (i.e. the
``confidence level'') is a tiny $1.965 \times 10^{-355}$.  Moreover,
the ratio of the last two recovered integer coefficients is precisely
$-11!  = -39916800$.  Given these facts, we can dismiss the
possibility that the recovered relation is a spurious numerical
artifact.

\vspace{1ex} \noindent {\bf Bifurcation to a 16-cycle:}
A second large application that we shall mention here is the problem
of determining the polynomial satisfied by the constant $B_4 =
3.564407268705 \cdots$, the fourth bifurcation point of the logistic
map $x_{k+1} = r x_k (1 - x_k)$. In section 3 we noted that an 8-cycle
begins at $r = B_3$, where $B_3$ satisfies a polynomial equation of
degree 12.  At $r = B_4$, this gives way to 16-cycle.  It has been
recognized that all $B_k$ are algebraic, but nothing has been known
about the degrees or the coefficients of the polynomials satisfied by
these constants for $k > 3$.  Some conjectural reasoning had suggested
that $B_4$ might satisfy a 240-degree polynomial, and some further
analysis had suggested that the constant $\alpha = - B_4 (B_4 - 2)$
might satisfy a 120-degree polynomial.  In order to test this
hypothesis, the three-level multi-pair program was applied to the
121-long vector $(1, \alpha, \alpha^2, \cdots, \alpha^{120})$.

In this case the input data was computed to over 10,000 digit
accuracy, and a working precision of 9,500 digits was employed in the
three-level multi-pair program.  A relation was detected at iteration
56,666.  The minimum and maximum $y$ vector entries at the point of
detection were $1.086 \times 10^{-9428}$ and $3.931 \times
10^{-8889}$, which form the ratio $2.763 \times 10^{-540}$.  Further,
the recovered integer coefficients descend monotonically from $257^{30}
\approx 1.986 \times 10^{72}$ to one.  Again, these facts argue very
strongly against the solution being a spurious numerical artifact.

\vspace{1ex} \noindent {\bf Reductions to Multiple Clausen Values:}
As a third application, consider sums of the form
\begin{eqnarray*}
S(k)&:=&\sum_{n>0}\frac{1}{n^{k}{2n \choose n}}\label{S}
\end{eqnarray*}
with, for example, $S(4)=17\pi^4/3240$.
Researchers have sought analytic evaluations of these constants for $k > 4$.
As a result of PSLQ computations, the constants
$\{S(k)\mid k=5\ldots20\}$ have been evaluated in terms of multiple
zeta values and {\em multiple Clausen values\/} of the form \cite{BorwBroad}
\begin{eqnarray*}
M(a,b) &:=& \sum_{n_1>n_2>\ldots>n_b>0}
  \frac{\sin(n_1 \pi /3)}{n_1^a} \prod_{j=1}^b\frac{1}{n_j}
\end{eqnarray*}
with, for example,
\begin{eqnarray*}
S(9) &=& \pi\left[2M(7,1) + \frac83M(5,3) + \frac89\zeta(2)M(5,1)\right]
  - \frac{13921}{216}\zeta(9) \\
  &&{} + \frac{6211}{486}\zeta(7)\zeta(2)
  + \frac{8101}{648}\zeta(6)\zeta(3)
  + \frac{331}{18}\zeta(5)\zeta(4) - \frac89\zeta^3(3)
\end{eqnarray*}

\noindent
The evaluation of the constant $S(20)$ is a 118-dimensional integer
relation problem, which required 4800 digit arithmetic.  In this case
a relation was detected at iteration 27,531.  The minimum and maximum
$y$ vector entry at detection were $7.170 \times 10^{-4755}$ and
$3.513 \times 10^{-4375}$, which gives a confidence ratio of $2.040
\times 10^{-380}$.  The actual solution for this problem is shown in
Table 5.  In this table, irreducible multiple zeta values such
$\zeta(5,3):=\sum_{j>k>0} j^{-5} k^{-3}$ occur.  Moreover, there are
alternating Euler sums, such as
$\zeta(\overline9,\overline3):=\sum_{j>k>0} (-1)^j j^{-9} \,(-1)^k
k^{-3}$, where an alternating sign is indicated by a bar.  The
presence of the latter results from another discovery obtained with
PSLQ \cite{BGK}, namely that some multiple zeta values may be reduced
to alternating Euler sums with fewer summations.  Finally, the
combinations \cite{BorwBroad}
\[\zeta_A(a,b,c):=
\zeta(\overline{a},\overline{b},c)
+\zeta(\overline{a},b,\overline{c})
+\zeta(a,\overline{b},\overline{c})
\]
serve to reduce 5-fold multiple zeta values to 3-fold alternating Euler
sums.

\begin{table}
\begin{center}
\begin{tiny}
$$\begin{array}{|l|}\hline
 525990827847624469523748125835264000
\,S(20) = \\[3pt]{}
-15024402006639545347476341466358480896000
\,\pi\,M(17,2)
+614357286926025380403737810975588352000
\,\pi\,M(15,4)\\[2pt]{}
-33663412982247966049519880053456896000
\,\pi\,M(13,6)
+204785762308675126801245936991862784000
\,\pi\,M(15,2)\,\zeta(2)\\[2pt]{}
-11221137660749322016506626684485632000
\,\pi\,M(13,4)\,\zeta(2)
-7792456708853695844796268530892800000
\,\pi\,M(13,2)\,\zeta(4)\\[2pt]{}
+65832426829545801661197345390290033253800417
\,\zeta(20)
-1655150248639886171642409815524246277640960
\,\zeta(17,3)\\[2pt]{}
-87407857867972646063318792204545819545600
\,\zeta(17)\,\zeta(3)
+239001490518032437117759318070284363571904
\,\zeta(15,5)\\[2pt]{}
+6475497072134876357497140759587182503936
\,\zeta(15,3)\,\zeta(2)
-11343388910891633971745524946475581811513600
\,\zeta(15)\,\zeta(5)\\[2pt]{}
+76505310594054968968541596301477435326464
\,\zeta(15)\,\zeta(3)\,\zeta(2)
-5427506872793330621343984298741119861120
\,\zeta(14)\,\zeta^2(3)\\[2pt]{}
-33725186900885181072542216636542494977280
\,\zeta(13,7)
-1079236594149043072329862323338197518336
\,\zeta(13,5)\,\zeta(2)\\[2pt]{}
-50485931801186079342010895425290633062400
\,\zeta(13,3)\,\zeta(4)
-24430610879956273104022963748303711510447040
\,\zeta(13)\,\zeta(7)\\[2pt]{}
+796530831594947602965411064203762718396416
\,\zeta(13)\,\zeta(5)\,\zeta(2)
-48476322702940293939397763722185147340800
\,\zeta(13)\,\zeta(4)\,\zeta(3)\\[2pt]{}
-2459446142542578280833853163647795200
\,\zeta(12)\,\zeta(5,3)
-7183917419981873615355846546110107008000
\,\zeta(12)\,\zeta(5)\,\zeta(3)\\[2pt]{}
+6554036738326690659991123688262156748800
\,\zeta(11,5)\,\zeta(4)
-674581129238392279111385274785342054400
\,\zeta(11,3,3,3)\\[2pt]{}
+155743130140661296228413518954716262400
\,\zeta(11,3,3)\,\zeta(3)
+13856996845301527891423305382301558784000
\,\zeta(11,3)\,\zeta(6)\\[2pt]{}
+339959536740516778440799419126460108800
\,\zeta(11,3)\,\zeta^2(3)
-35543027806069609369237745997797431835122560
\,\zeta(11)\,\zeta(9)\\[2pt]{}
+1912599458053045671374932296869893271531520
\,\zeta(11)\,\zeta(7)\,\zeta(2)
-8624509220693012537969847600322793702400
\,\zeta(11)\,\zeta(6)\,\zeta(3)\\[2pt]{}
-159424648200337153322748394462255349760000
\,\zeta(11)\,\zeta(5)\,\zeta(4)
-386372041666595966843560058208603955200
\,\zeta(11)\,\zeta^3(3)\\[2pt]{}
-4526144521471219675886040639917260800
\,\zeta(10)\,\zeta(7,3)
-4235684121072319605030836248657970626560
\,\zeta(10)\,\zeta(7)\,\zeta(3)\\[2pt]{}
-4274427562442524135198151261132645652480
\,\zeta(10)\,\zeta^2(5)
+174910231480430088343102177690512998400
\,\zeta(9,5,3,3)\\[2pt]{}
+23201851844071266080584141499247820800
\,\zeta(9,5,3)\,\zeta(3)
-5388965272775430297200443154254448394240
\,\zeta(9,5)\,\zeta(6)\\[2pt]{}
-96144480802344282256962346694615654400
\,\zeta(9,5)\,\zeta^2(3)
-564799665543005814719751486159037931520
\,\zeta(9,3,5,3)\\[2pt]{}
+192405432086205157974414874044727296000
\,\zeta(9,3,3,3)\,\zeta(2)
-437636171132005416578131168531552665600
\,\zeta(9,3,3)\,\zeta(5)\\[2pt]{}
+82410232260928579141238186701396377600
\,\zeta(9,3,3)\,\zeta(3)\,\zeta(2)
-56820309831551194167334052913378508800
\,\zeta(9,3)\,\zeta(8)\\[2pt]{}
-78290750182007491017999587160883200
\,\zeta(9,3)\,\zeta(5,3)
+173223299338939829293781467642177536000
\,\zeta(9,3)\,\zeta(5)\,\zeta(3)\\[2pt]{}
-5389461879726322601463723747508224000
\,\zeta(9,3)\,\zeta^2(3)\,\zeta(2)
+1395360857314382550903663041050280719202304
\,\zeta^2(9)\,\zeta(2)\\[2pt]{}
-1543454230261900138881951172107169382400
\,\zeta(9)\,\zeta(8)\,\zeta(3)
-511939532590839950285975762448130830336000
\,\zeta(9)\,\zeta(7)\,\zeta(4)\\[2pt]{}
+89785104680812821069278191239404195328000
\,\zeta(9)\,\zeta(6)\,\zeta(5)
-1309132727087420901925773113189990400
\,\zeta(9)\,\zeta(5,3)\,\zeta(3)\\[2pt]{}
+1731994708600523066371520212640192716800
\,\zeta(9)\,\zeta(5)\,\zeta^2(3)
-1309132727087420901925773113189990400
\,\zeta(9)\,\zeta(3,5,3)\\[2pt]{}
+49863866344508636305947249931911168000
\,\zeta(9)\,\zeta^3(3)\,\zeta(2)
-13797482183512283560940162429818580121600
\,\zeta(8)\,\zeta(7)\,\zeta(5)\\[2pt]{}
+20525248296522064841059215485763379200
\,\zeta(8)\,\zeta^4(3)
+533245759266957435712480647773027635200
\,\zeta(7,7,3,3)\\[2pt]{}
+39157503832984121716572521716488652800
\,\zeta(7,7,3)\,\zeta(3)
+223377519430349618539918571265503416320
\,\zeta(7,5,5,3)\\[2pt]{}
-23700768289019448234404348103552000000
\,\zeta(7,5,5)\,\zeta(3)
-184392479550115407127835175133101465600
\,\zeta(7,5,3,5)\\[2pt]{}
-192646077208687087875906081369587712000
\,\zeta(7,5,3,3)\,\zeta(2)
+609805989326475901096307023020578611200
\,\zeta(7,5,3)\,\zeta(5)\\[2pt]{}
-92883012339775157113672910718900633600
\,\zeta(7,5,3)\,\zeta(3)\,\zeta(2)
-1083004232781819170351004903486259200
\,\zeta(7,3,5,3)\,\zeta(2)\\[2pt]{}
+642228810086780199757027863429120000
\,\zeta(7,3,3)\,\zeta(7)
+1715165074541342577478967041720320000
\,\zeta(7,3,3)\,\zeta(5)\,\zeta(2)\\[2pt]{}
-256615593289239779065106729533440000
\,\zeta^2(7,3)
-117611319397120633884272025717964800
\,\zeta(7,3)\,\zeta(5,3)\,\zeta(2)\\[2pt]{}
+266618082812610684431891803668480000
\,\zeta(7,3)\,\zeta^2(5)
+68054784737327925282389654900302533120000
\,\zeta^2(7)\,\zeta(6)\\[2pt]{}
+1571090732393362601235892759430587238400
\,\zeta^2(7)\,\zeta^2(3)
-1016434130097344129482122765187153920
\,\zeta(7)\,\zeta(5,5,3)\\[2pt]{}
-293590313182528091317498451853312000
\,\zeta(7)\,\zeta(5,3)\,\zeta(5)
+1133345755060987206065174842254330470400
\,\zeta(7)\,\zeta^2(5)\,\zeta(3)\\[2pt]{}
+62419983317149231400720825830146048000
\,\zeta(7)\,\zeta(5)\,\zeta^2(3)\,\zeta(2)
+101029201230288166627783621503025152000
\,\zeta(7)\,\zeta(4)\,\zeta^3(3)\\[2pt]{}
-3927398181262262705777319339569971200
\,\zeta(6)\,\zeta(5,3,3,3)
+88122944806249232884839801748733952000
\,\zeta(6)\,\zeta(5)\,\zeta^3(3)\\[2pt]{}
-1963699090631131352888659669784985600
\,\zeta(6)\,\zeta(3,5,3)\,\zeta(3)
+1239443914180202003982888789718597632
\,\zeta(5,5,5,3)\,\zeta(2)\\[2pt]{}
-2714534591307431519290045171362693120
\,\zeta(5,5,3)\,\zeta(5)\,\zeta(2)
-7854796362524525411554638679139942400
\,\zeta(5,3,3,3,3,3)\\[2pt]{}
+7854796362524525411554638679139942400
\,\zeta(5,3,3,3,3)\,\zeta(3)
-3927398181262262705777319339569971200
\,\zeta(5,3,3,3)\,\zeta^2(3)\\[2pt]{}
-117611319397120633884272025717964800
\,\zeta(5,3)\,\zeta^2(5)\,\zeta(2)
+327283181771855225481443278297497600
\,\zeta(5,3)\,\zeta^4(3)\\[2pt]{}
-310534753804603441554226729609432166400
\,\zeta^4(5)
+46595661441120443976355182952120320000
\,\zeta^3(5)\,\zeta(3)\,\zeta(2)\\[2pt]{}
+151543801845432249941675432254537728000
\,\zeta^2(5)\,\zeta(4)\,\zeta^2(3)
-888340064809321326306774612521779200
\,\zeta(5)\,\zeta^5(3)\\[2pt]{}
-654566363543710450962886556594995200
\,\zeta(3,5,3)\,\zeta^3(3)
+15584913417707391689592537061785600
\,\zeta^6(3)\,\zeta(2)\\[2pt]{}
+31338860750207444579474396657221632000
\,\zeta(8)\,\zeta(\overline9,\overline3)
+18542546095738616293736744327577600000
\,\zeta(5)\,\zeta(3)\,\zeta(\overline9,\overline3)\\[2pt]{}
+8537710593272460830117524829896704000
\,\zeta^2(3)\,\zeta(2)\,\zeta(\overline9,\overline3)
+675871149225360968655980361248931840000
\,\zeta(4)\,\zeta(\overline{13},\overline3)\\[2pt]{}
+254015007537749154775389871602030084096
\,\zeta(2)\,\zeta(\overline{15},\overline3)
-1692980876937872291412185599949615923200
\,\zeta(\overline{17},\overline3)\\[2pt]{}
-12361697397159077529157829551718400000
\,\zeta(5)\,\zeta_A(9,3,3)
-11383614124363281106823366439862272000
\,\zeta(3)\,\zeta(2)\,\zeta_A(9,3,3)\\[2pt]{}
+212786017863098254535236772683579392000
\,\zeta(3)\,\zeta_A(11,3,3)
+174238991699437976124847445508096000
\,\zeta(2)\,\zeta(6,\overline5,\overline4,3)\\[2pt]{}
+65242291818339575848332989300736000
\,\zeta(8,\overline5,\overline4,3)
+103014144976325646076315246264320000
\,\zeta(6,\overline5,\overline6,3)\\[2pt]\hline
\end{array}$$
\end{tiny}
\end{center}
\caption{Solution for $S(20)$ found with the three-level program}
\end{table}

\newpage

\vspace{1ex} \noindent {\bf Parallelization:}
These three problems were first solved by the second author running a
three-level implementation of PSLQ on a DecAlpha machine at the Open
University, with a single 433 MHz processor, and 1 Gbyte of main
memory.  They were then used as benchmarks for a multiprocessor
version of the new three-level multi-pair program, using the OpenMP
programming model, on a 64-CPU SGI Origin-2000 system at the Lawrence
Berkeley Laboratory.  Run times are given in Table 6.  Timings on 48
processors show a speedup of 19.40 times on the Fibonacci conjecture
problem, 22.44 times on the $B_4$ problem, and 17.81 times on the
$S(20)$ problem.  Given the challenge of very limited concurrency
inherent in this type of calculation, we are encouraged by these
figures.

\begin{table}[h]
\begin{center}
\begin{tabular}{|l|r|r|r|r|r|r|r|}
\hline
& \multicolumn{2}{|c|}{Fibonacci} & \multicolumn{2}{|c|}{$B_4$} &
  \multicolumn{2}{|c|}{$S(20)$} \\
Processors & Time & Speedup & Time & Speedup & Time & Speedup \\
\hline
 1 & 47788 &  1.00 & 90855 &  1.00 & 23208 &  1.00 \\
 2 & 24665 &  1.94 & 46134 &  1.97 & 11973 &  1.94 \\
 4 & 12945 &  3.69 & 23966 &  3.79 &  6305 &  3.68 \\
 8 &  7076 &  6.75 & 12924 &  7.03 &  3470 &  6.69 \\
16 &  4180 & 11.43 &  7424 & 12.24 &  2126 & 10.92 \\
32 &  2994 & 15.96 &  4865 & 18.68 &  1548 & 14.99 \\
48 &  2463 & 19.40 &  4049 & 22.44 &  1303 & 17.81 \\
\hline
\end{tabular} \end{center}
\caption{Timings for three large problems using the parallel
three-level multi-pair program}
\end{table}

\vspace{2ex} \noindent
{\bf 9. Conclusion}

We have accelerated the conventional implementation of the PSLQ
algorithm in three ways.  First, we utilized a two-level and a
three-level scheme, which permit most if not all iterations to be
performed using ordinary 64-bit double precision arithmetic, and
updating the multiprecision arrays only as needed.  This resulted in a
speedup of up to 65 times over the straightforward one-level program.
Secondly, we developed a new integer relation algorithm, a variant of
PSLQ that we have termed the ``multi-pair'' algorithm.  We also
demonstrated two-level and three-level implementations of this new
algorithm.  These techniques resulted in an additional speedup of up
to 22\%, comparing the three-level multi-pair program to the
three-level PSLQ program.  Finally, we showed how that this new
algorithm, unlike PSLQ, is reasonably well suited for parallel
processing.  We demonstrated a parallel three-level implementation of
the multi-pair algorithm that achieved an {\em additional\/} speedup
of up to 22 times.

In consequence, we are able to discover in days relations that would
previously have taken years to unveil.  We applied the programs to
three large problems, obtaining results not previously known in the
literature.  We believe that these demonstrations open up a novel way
of doing science.  We are confident that many more discoveries can be
made in this manner.

\raggedright

\end{document}